\documentclass[reqno,12pt]{amsart}
\usepackage{amscd,amssymb}

\numberwithin{equation}{section}
\theoremstyle{plain}

\theoremstyle{plain}
\newtheorem{Thm}[subsection]{Theorem}
\newtheorem{Cor}[subsection]{Corollary}
\newtheorem{Lem}[subsection]{Lemma}
\newtheorem{Prop}[subsection]{Proposition}

\theoremstyle{definition}
\newtheorem{Def}[subsection]{Definition}

\theoremstyle{remark}

\newtheorem{Rem}[subsection]{Remark}

\newif\ifShowLabels
\ShowLabelstrue
\newdimen\theight
\def\TeXref#1{%
        \leavevmode\vadjust{\setbox0=\hbox{{\tt
                \quad\quad  {\small \textrm #1}}}%
        \theight=\ht0
        \advance\theight by \lineskip
        \kern -\theight \vbox to
        \theight{\rightline{\rlap{\box0}}%
        \vss}%
        }}%

\ShowLabelsfalse

\renewcommand{\sec}[2]{\section{#2}\label{S:#1}%
        \ifShowLabels \TeXref{{S:#1}} \fi}
\newcommand{\ssec}[2]{\subsection{#2}\label{SS:#1}%
        \ifShowLabels \TeXref{{SS:#1}} \fi}

\newcommand{\refs}[1]{Section ~\ref{S:#1}}
\newcommand{\refss}[1]{Subsection ~\ref{SS:#1}}
\newcommand{\reft}[1]{Theorem ~\ref{T:#1}}
\newcommand{\refl}[1]{Lemma ~\ref{L:#1}}
\newcommand{\refp}[1]{Proposition ~\ref{P:#1}}
\newcommand{\refc}[1]{Corollary ~\ref{C:#1}}

\newcommand{\refe}[1]{\eqref{E:#1}}

\newenvironment{thm}[1]%
        { \begin{Thm} \label{T:#1}  \ifShowLabels \TeXref{T:#1} \fi }%
        { \end{Thm} }

\renewcommand{\th}[1]{\begin{thm}{#1} \sl }
\renewcommand{\eth}{\end{thm} }

\newenvironment{lemma}[1]%
        { \begin{Lem} \label{L:#1}  \ifShowLabels \TeXref{L:#1} \fi }%
        { \end{Lem} }
\newcommand{\lem}[1]{\begin{lemma}{#1} \sl}
\newcommand{\elem}{\end{lemma}}

\newenvironment{propos}[1]%
        { \begin{Prop} \label{P:#1}  \ifShowLabels \TeXref{P:#1} \fi }%
        { \end{Prop} }
\newcommand{\prop}[1]{\begin{propos}{#1}\sl }
\newcommand{\eprop}{\end{propos}}

\newenvironment{corol}[1]%
        { \begin{Cor} \label{C:#1}  \ifShowLabels \TeXref{C:#1} \fi }%
        { \end{Cor} }
\newcommand{\cor}[1]{\begin{corol}{#1} \sl }
\newcommand{\ecor}{\end{corol}}

\newenvironment{defeni}[1]%
        { \begin{Def} \label{D:#1}  \ifShowLabels \TeXref{D:#1} \fi }%
        { \end{Def} }
\newcommand{\defe}[1]{\begin{defeni}{#1} \sl }
\newcommand{\edefe}{\end{defeni}}

\newenvironment{remark}[1]%
        { \begin{Rem} \label{R:#1}  \ifShowLabels \TeXref{R:#1} \fi }%
        { \end{Rem} }
\newcommand{\rem}[1]{\begin{remark}{#1}}
\newcommand{\erem}{\end{remark}}

\newcommand{\eq}[1]%
        { \ifShowLabels \TeXref{E:#1} \fi
           \begin{equation} \label{E:#1} }
\newcommand{\eeq}{\end{equation}}

\newcommand{\prf}{ \begin{proof} }
\newcommand{\eprf}{ \end{proof} }

\newcommand\alp{\alpha}         
\newcommand\bet{\beta}
\newcommand\gam{\gamma}         \newcommand\Gam{\Gamma}
\newcommand\del{\delta}         \newcommand\Del{\Delta}
\newcommand\eps{\varepsilon}

\newcommand\iot{\iota}

\newcommand\lam{\lambda}                \newcommand\Lam{\Lambda}
\newcommand\sig{\sigma}

\newcommand\calA{{\mathcal{A}}}

\newcommand\calC{{\mathcal{C}}}

\newcommand\calE{{\mathcal{E}}}
\newcommand\calF{{\mathcal{F}}}

\newcommand\calL{{\mathcal{L}}}

\newcommand\calO{{\mathcal{O}}}

\newcommand\calS{{\mathcal{S}}}

\newcommand\calW{{\mathcal{W}}}

\newcommand\calZ{{\mathcal{Z}}}

\newcommand\bfc{{\mathbf c}}

\newcommand\ZZ{\mathbb{Z}}

\newcommand\CC{\mathbb{C}}

\newcommand\nek{,\ldots,}
\newcommand\sdp{\times \hskip -0.3em {\raise 0.3ex
\hbox{$\scriptscriptstyle |$}}} 

\newcommand\Dom{\operatorname{Dom}}

\newcommand\End{\operatorname{End\,}}

\newcommand\Hom{\operatorname {Hom}}

\newcommand\ind{\operatorname{ind}}

\newcommand\Ker{\operatorname{Ker}}

\newcommand\RE{\operatorname{Re}}

\newcommand{\sign}{\operatorname{sign}}

\newcommand\spin{\operatorname{spin}}

\newcommand\Tr{\operatorname{Tr}}

\renewcommand\oe{{\overline{e}}}

\newcommand\oz{{\overline{z}}}

\newcommand\tilD{{\widetilde{D}}}

\newcommand\tilJ{{\widetilde{J}}}

\newcommand\tilS{{\widetilde{S}}}

\newcommand\tilDel{{\widetilde{\Delta}}}

\newcommand\tileps{{\widetilde{\eps}}}

\renewcommand{\>}{\rangle}
\newcommand{\<}{\langle}

\renewcommand{\d}{\text{\( \partial\)}}
\newcommand{\p}{\bar{\d}}

\newcommand{\n}{\nabla}
\newcommand{\tD}{\tilDel}
\newcommand{\tJ}{\tilJ}

\newcommand{\nk}{\nabla^{\E\otimes\L^k}}
\newcommand{\Ek}{{\E\otimes\L^k}}
\newcommand{\tilnk}{\widetilde{\nabla}^{\E\otimes\L^k}}
\newcommand{\tildk}{\widetilde{\Del}_k}
\newcommand{\E}{\calE}\newcommand{\W}{\calW}
\newcommand{\Z}{\calZ}\newcommand{\F}{\calF}
\newcommand{\A}{\calA}\renewcommand{\O}{\calO}
\renewcommand{\L}{\calL}
\renewcommand{\tilS}{\widetilde{\calS}}
\newcommand{\Lk}{\L^k}
\newcommand{\fes}{F^{\E/\calS}}

\newcommand{\EndC}{\operatorname{End}_{C(M)}\,}
\newcommand{\g}{{\Gam}}
\newcommand{\gc}{{\Gam(M,C(M))}}
\newcommand{\gme}{{\Gam(M,\E)}}
\newcommand{\ha}{^{1,0}}
\newcommand{\ah}{^{0,1}}

\newcommand{\hor}{^{\text{hor}}}
\renewcommand{\vert}{^{\text{vert}}}

\newcommand{\ka}{K\"ahler }
\begin{document}

\title{Vanishing theorems for the half-kernel of a Dirac operator}
\author{Maxim Braverman}
\address{Institute of Mathematics\\
        The Hebrew University   \\
         Jerusalem 91904 \\
         Israel
         }
\email{maxim@math.huji.edu}
\thanks{This research was partially supported by grant No. 96-00210/1 from
the United States-Israel Binational Science Foundation (BSF)}

\begin{abstract}
We obtain a vanishing theorem for the half-kernel of a Dirac operator on a Clifford
module twisted by a sufficiently large power of a line bundle, whose curvature
is non-degenerate at any point of the base manifold. In particular, if the base
manifold is almost complex, we prove a vanishing theorem for the half-kernel of a
$\spin^c$ Dirac operator twisted by a line bundle with curvature of a mixed
sign. In this case we also relax the assumption of non-degeneracy of the
curvature. These results are generalization of a vanishing theorem of Borthwick
and Uribe. As an application we obtain a new proof of the classical
Andreotti-Grauert vanishing theorem for the cohomology of a compact complex
manifold with values in the sheaf of holomorphic sections of a holomorphic
vector bundle, twisted by a large power of a holomorphic line bundle with
curvature of a mixed sign.

As another application we calculate the sign of the index of a signature operator
twisted by a large power of a line bundle.
\end{abstract}

\maketitle
\tableofcontents

\sec{introd}{Introduction}

One of the most fundamental facts of complex geometry is the Kodaira vanishing
theorem for the cohomology of the sheaf of sections of a holomorphic vector
bundle twisted by a large power of a positive line bundle. In 1962,
Andreotti and Grauert \cite{AndGr62} obtained the following
generalization of this result to the case when the line bundle is not
necessarily positive. Let $\L$ be a holomorphic line bundle over a compact
complex $n$-dimensional manifold $M$. Suppose $\L$ admits a holomorphic
connection whose curvature $F^\L$ has a least $q$ negative and at least $p$
positive eigenvalues at any point of $M$. Then the Andreotti-Grauert theorem
asserts that, for any holomorphic vector bundle $\W$ over $M$, the cohomology
$H^j(M,\O(\W\otimes\Lk))$ of $M$ with coefficients in the sheaf of holomorphic
sections of the tensor product $\W\otimes\Lk$ vanishes for $k\gg0,
\ j\not=q,q+1\nek n-p$. In particular, if $F^\L$ is non-degenerate at all points
of $M$, then the number $q$ of negative eigenvalues of $F^\L$ is independent of
$x\in M$, and the Andreotti-Grauert theorem implies that the cohomology
$H^j(M,\O(\W\otimes\Lk))$ vanishes for $k\gg0, j\not=q$.

If $M, \W$ and $\L$ are endowed with metrics, then the cohomology
$H^*(M,\O(\W\otimes\Lk))$ is isomorphic to the kernel of the Dolbeault-Dirac
operator
\[
        \sqrt2(\p+\p^*):\A^{0,*}(M,\W\otimes\Lk)\to
                                \A^{0,*}(M,\W\otimes\Lk).
\]
Here $\A^{0,*}(M,\W\otimes\Lk)$ denotes the space of $(0,*)$-differential forms
on $M$ with values in $\W\otimes\Lk$. The Andreotti-Grauert theorem implies, in
particular, that the restriction of the kernel of the Dolbeault-Dirac operator
on the space $\A^{0,\text{odd}}(M,\W\otimes\Lk)$ (resp.
$\A^{0,\text{even}}(M,\W\otimes\Lk)$) vanishes provided the curvature $F^\L$ is
non-degenerate and has an even (resp. an odd) number of negative eigenvalues at
any point of $M$.

The last statement may be extended to the case when the manifold $M$ is not
complex. First step in this direction was done by Borthwick and Uribe
\cite{BoUr96}, who showed that, if $M$ is an almost \ka manifold and $\L$ is a
positive line bundle over $M$, then the restriction of the kernel of the
$\spin^c$-Dirac operator $D_k:\A^{0,*}(M,\Lk)\to\A^{0,*}(M,\Lk)$ on the space
$\A^{0,\text{odd}}(M,\W\otimes\Lk)$ vanishes for $k\gg0$. Moreover, they showed
that, for any $\alp\in \Ker D_k$, ``most of the norm" of $\alp$ is concentrated
in $\A^{0,0}(M,\Lk)$. This result generalizes the Kodaira vanishing theorem to
the case of an almost \ka manifolds.

One of the results of the present paper is the extension of the Borthwick-Uribe
theorem to the case when the curvature $F^\L$ of $\L$ is not positive. In other
words, we extend the Andreotti-Grauert theorem to almost complex manifolds.

More generally, assume that $M$ is a compact oriented even-dimensional
Riemannian manifold and let $C(M)$ denote the Clifford bundle of $M$, i.e., a
vector bundle whose fiber at any point is isomorphic to the Clifford algebra of
the cotangent space. Let $\E$ be a self-adjoint Clifford module over $M$, i.e.,
a Hermitian vector bundle over $M$ endowed with a fiberwise action of $C(M)$.
Then (cf. \refss{chir}) $\E$ possesses a natural grading $\E=\E^+\oplus\E^-$.
Let $\L$ be a Hermitian line bundle endowed with a Hermitian connection $\n^\L$
and let $\E$ be a Hermitian vector bundle over $M$ endowed with an Hermitian
connection $\n^\E$. These data define (cf. \refs{dirac}) a self-adjoint Dirac
operator $D_k:\Gam(\E\otimes\Lk)\to
\Gam(\E\otimes\Lk)$. The curvature  $F^\L$  of  $\n^\L$  is an imaginary valued
2-form on $M$. If it is non-degenerate at all points of $M$, then $iF^\L$ is a
symplectic form on $M$, and, hence, defines an orientation of $M$. Our main
result (\reft{main}) states that {\em the restriction of the kernel of $D_k$ to
$\Gam(\E^-\otimes\Lk)$ (resp. to $\Gam(\E^+\otimes\Lk))$ vanishes for large $k$
if this orientation coincides with (resp. is opposite to) the given orientation of
$M$.}

Our result may be considerably refined when $M$ is an almost complex
$2n$-dimensional manifold and the curvature $F^\L$ is a $(1,1)$-form on $M$. In
this case, $F^\L$ may be considered as a sesquilinear form on the holomorphic
tangent bundle to $M$. Let $\W$ be a Hermitian vector bundle over $M$ endowed
with an Hermitian connection. Then (cf. \refss{almcomp}) there is a canonically
defined Dirac operator $D_k:\A^{0,*}(M,\W\otimes\Lk)\to
\A^{0,*}(M,\W\otimes\Lk)$. We prove (\reft{UB+}) that {\em if $F^\L$ has at
least $q$ positive and at least $p$ negative eigenvalues at every point of $M$,
then, for large $k$, ``most of the norm" of any element $\alp\in \Ker D_k$ is
concentrated in $\bigoplus_{j=q}^{n-p}\A^{0,j}(M,\W\otimes\Lk)$.} In particular,
{\em if the sesquilinear form $F^\L$ is non-degenerate and has exactly $q$
negative eigenvalues at any point of $M$, then ``most of the norm" of $\alp\in
\Ker D_k$ is concentrated in $\A^{0,q}(M,\W\otimes\Lk)$, and, depending on the
parity of $q$, the restriction of the kernel of $D_k$ either to
$\A^{0,\text{odd}}(M,\W\otimes\Lk)$ or to $\A^{0,\text{even}}(M,\W\otimes\Lk)$
vanishes.} These results generalize both the Andreotti-Grauert and the
Borthwick-Uribe vanishing theorems. In particular, we obtain a new proof of the
Andreotti-Grauert theorem.

As another application of \reft{main}, we study the index of a signature
operator twisted by a line bundle having a non-degenerate curvature. We prove
(\refc{derham}) that, {\em if the orientation defined by the curvature of $\L$
coincides with (resp. is opposite to) the given orientation of $M$, then this
index is non-negative (resp. non-positive).}

The proof of our main vanishing theorem (\reft{main}) is based on an estimate of
the square $D_k^2$ of the twisted Dirac operator for large values of $k$. This
estimate is obtained in two steps. First we use the Lichnerowicz formula to
compare $D_k^2$ with the metric Laplacian
$\Del_k=\left(\n^{\E\otimes\Lk}\right)^*\n^{\E\otimes\Lk}$. Then we use the
method of \cite{GuiUr88,BoUr96} to estimate the large $k$ behavior of the metric
Laplacian.

\subsection*{Contents}
The paper is organized as follows:

In \refs{dirac}, we briefly recall some basic facts about Clifford modules and Dirac
operators. We also present some examples of Clifford modules which are used in the rest
of the paper.

In \refs{main}, we formulate the main results of the paper and discuss their
applications. The rest of the paper is devoted to the proof of these results.

In \refs{prmain}, we present the proof of \reft{main} (the vanishing theorem for
the half-kernel of a Dirac operator). The proof is based on two statements
(Propositions~\ref{P:D-D} and \ref{P:lapl}) which are proven in the later
sections.

In \refs{prUB}, we prove an estimate on the Dirac operator on an almost complex
manifold (\refp{estDk}) and use it to prove \reft{UB+} (our analogue of the
Andreotti-Grauert vanishing theorem for almost complex manifolds). The proof is
based on Propositions~\ref{P:lapl} and \ref{P:D-Dac} which are proved in later
sections.

In \refs{prAG}, we prove the Andreotti-Grauert theorem (\reft{AG}).

In \refs{lichn}, we use the Lichnerowicz formula to prove Propositions~\ref{P:D-D},
\ref{P:D-Dac} and \ref{P:gr}. These results establish the connection between the Dirac
operator and the metric Laplacian. They are used in the proofs of
Theorems~\ref{T:main}, \ref{T:UB+} and \ref{T:AG}.

{}Finally, in \refs{lapl}, we apply the method of \cite{GuiUr88,BoUr96} to prove
\refp{lapl} (the estimate on the metric Laplacian).

\subsection*{Acknowledgments}
This paper grew out of our joint effort with Yael Karshon to understand the
possible applications and generalizations of the results of \cite{BoUr96}. I
would like to thank Yael Karshon for drawing my attention to this paper and for
her help and support during the work on this project. Many of the results of the
present paper were obtained in our discussions with Yael.

I would like to thank Joseph Bernstein and Gideon Maschler for valuable
discussions.

\sec{dirac}{Clifford modules and Dirac operators}

In the first part of this section we briefly recall the definitions and some basic
facts about Clifford modules and Dirac operators. We refer the reader to
\cite{BeGeVe,Duis96,LawMic89} for details. In our exposition we adopt the notations of
\cite{BeGeVe}.

In the second part of the section we present some examples of Clifford modules, which
will be used in the subsequent sections.

\ssec{clbundle}{The Clifford bundle}
Suppose $(M,g)$ is an oriented Riemannian manifold of dimension $2n$. For any
$x\in M$, we denote by $C(T^*_xM)=C^+(T^*_xM)\oplus C^-(T^*_xM)$ the Clifford
algebra of the cotangent space $T_x^*M$, cf. \cite[\S3.3]{BeGeVe}.

The {\em Clifford bundle \/} $C(M)$ of $M$ (cf. \cite[\S3.3]{BeGeVe}) is the
$\ZZ_2$-graded bundle over $M$, whose fiber at $x\in M$ is $C(T^*_xM)$.

The Riemannian metric $g$ induces the Levi-Civita connection $\n$ on $C(M)$
which is compatible with the multiplication and preserves the $\ZZ_2$-grading on
$C(M)$.

\ssec{clmodule}{Clifford modules}
A {\em Clifford module \/} on $M$ is a complex vector bundle $\E$ on $M$ endowed
with an action of the bundle $C(M)$. We write this action as
\[
        (a,s) \ \mapsto \ c(a)s, \quad \mbox{where} \quad
                        a\in \gc, \ s\in \gme.
\]

A Clifford module $\E$ is called {\em self-adjoint \/} if it is endowed with a
Hermitian metric such that the operator $c(v):\E_x\to\E_x$ is skew-adjoint, for
any $x\in M$ and any $v\in T_x^*M$.

A connection $\n^\E$ on a Clifford module $\E$ is called a {\em Clifford
connection \/} if it is {\em compatible with the Clifford action}, i.e., if for
any $a\in \gc$ and $X\in\g(M,TM)$,
\[
        [\n^\E_X,c(a)] \ = \ c(\n_X a).
\]
In this formula, $\n_X$ is the Levi-Civita covariant derivative on $C(M)$, and
$[\n^\E_X,c(a)]$ denotes the commutator of the operators $\n_X^\E$ and $c(a)$.

Suppose $\E$ is a Clifford module and $\calW$ is a vector bundle over $M$. The
{\em twisted Clifford module obtained from $\E$ by twisting with $\calW$ \/} is
the bundle $\E\otimes\calW$ with Clifford action $c(a)\otimes1$. Note that the
twisted Clifford module $\E\otimes\calW$ is self-adjoint if and only if so is
$\E$.

Let $\n^\calW$ be a connection on $\calW$ and let $\n^\E$ be a Clifford
connection on $\E$. Then the {\em product connection}
\begin{equation}\label{E:nEW}
        \n^{\E\otimes\calW} \ = \ \n^\E\otimes 1 \ + \ 1\otimes \n^\calW
\end{equation}
is a Clifford connection on $\E\otimes\calW$.

\ssec{chir}{The chirality operator. The natural grading}
Let $e_1\nek e_{2n}$ be an oriented orthonormal basis of $C(T_x^*M)$ and
consider the element
\begin{equation}\label{E:Gam}
        \Gam \ = \ i^n\, e_1\cdots e_{2n} \ \in \ C(T_x^*M)\otimes\CC.
\end{equation}
This element is independent of the choice of the basis, anti-commutes with any
$v\in T_x^*M\subset C(T_x^*M)$, and satisfies $\Gam^2=-1$, cf.
\cite[\S3.2]{BeGeVe}. This element $\Gam$ is called the {\em chirality
operator}. We also denote by $\Gam$ the section of $C(M)$ whose restriction to
each fiber is equal to the chirality operator.

Let $\E$ be a Clifford module, i.e. (cf. \refss{clmodule}), a vector bundle over
$M$ endowed with a fiberwise action of $C(M)$. Set
\begin{equation}\label{E:grad}
        \E^\pm \ = \ \{v\in \E: \ \Gam\, v=\pm v\}.
\end{equation}
Then $\E=\E^+\oplus\E^-$ is a {\em graded module} over $C(M)$ in the sense that
$C^+(M)\cdot\E^\pm\subset \E^\pm$ and $C^-(M)\cdot\E^\pm\subset \E^\mp$.

We refer to the grading \refe{grad} as the {\em natural grading \/} on $\E$.
Note that this grading is preserved by any Clifford connection on $\E$. Also, if
$\E$ is a self-adjoint Clifford module (cf. \refss{clmodule}), then the
chirality operator $\Gam:\E\to\E$ is self-adjoint. Hence, the subbundles
$\E^\pm$ are orthogonal with respect to the Hermitian metric on $\E$. {\em In
this paper we endow all our Clifford modules with the natural grading.}

\ssec{dirac}{Dirac operators}
The {\em Dirac operator \/} associated to a Clifford connection $\n^\E$ is
defined by the following composition
\begin{equation}\label{E:dir1}
\begin{CD}
        \gme @>\n^\E>> \g(M,T^*M\otimes \E) @>c>> \gme.
\end{CD}
\end{equation}
In local coordinates, this operator may be written as $D=\sum\,
c(dx^i)\,\n^\E_{\d_i}$. Note that $D$ sends even sections to odd sections and
vice versa: $D:\, \Gam(M,\E^\pm)\to \Gam(M,\E^\mp)$.

Suppose that the Clifford module $\E$ is endowed with a Hermitian structure and
consider the $L_2$-scalar product on the space of sections $\gme$ defined by the
Riemannian metric on $M$ and the Hermitian structure on $\E$. By
\cite[Proposition~3.44]{BeGeVe}, {\em the Dirac operator associated to a
Clifford connection $\n^\E$ is formally self-adjoint with respect to this scalar
product if and only if $\E$ is a self-adjoint Clifford module and $\n^\E$ is a
Hermitian connection.}

\

We finish this section with some examples of Clifford modules, which will be used
later.

\ssec{spin}{Spinor bundles}
Assume that $M$ is a spin manifold and let $\calS= \calS^+\oplus
\calS^-$  be a  {\em spinor bundle \/}  over  $M$ (cf. \cite[\S
3.3]{BeGeVe}). It is a minimal Clifford module in the sense that any other
Clifford module $\E$ may be decomposed as a tensor product
\begin{equation}\label{E:E=SW}
        \E \ =  \ \calS \ \otimes \calW,
\end{equation}
where $\calW= \Hom_{C(M)}(\calS,\E)$ and the action of the Clifford bundle
$C(M)$ is trivial on the second factor of \refe{E=SW}. In this case, the natural
grading $\E=\E^+\oplus\E^-$ is defined by $\E^\pm=\calS^\pm\otimes\W$.

The Riemannian metric on $M$ induces the Levi-Civita connection $\n^\calS$ on
$\calS$, which is compatible with the Clifford action. Moreover, a connection
$\n^\E$ on the twisted Clifford module $\E= \calS\otimes \calW$ is a Clifford
connection if and only if
\begin{equation}\label{E:nE=}
        \n^\E \ = \ \n^\calS\otimes 1 \ + \ 1\otimes \n^\calW
\end{equation}
for some connection $\n^\calW$ on $\calW$.

Note that locally the spinor bundle and, hence, the decompositions \refe{E=SW},
\refe{nE=} always exist. In particular, suppose that $\tilS$ is a Clifford
module whose fiber dimension is equal to $\dim\calS=2^n$. Then locally
$\tilS=\calS\otimes\calW$ for some locally defined complex line bundle $\calW$.
In this case $\tilS$ is called a $\spin^c$ vector bundle over $M$
(\cite[Ch.~5]{Duis96}, \cite[Appendix~D]{LawMic89}).

A Dirac operator on a $\spin^c$ vector bundle is called a {\em $\spin^c$ Dirac
operator}.

\ssec{difforms}{The exterior algebra}
Consider the exterior algebra $\Lam{T^*M}=\bigoplus_i\Lam^iT^*M$ of the
cotangent bundle $T^*M$. There is a canonical action of the Clifford bundle
$C(M)$ on $\Lam{T^*M}$ such that
\begin{equation}\label{E:cv}
        c(v)\, \alp \ = \ v\wedge\alp \ - \ \iot(v)\, \alp, \qquad
                        v\in \Gam(M,T^*M), \ \alp\in\Gam(M,\Lam{T^*M}).
\end{equation}
Here $\iot(v)$ denotes the contraction with the vector $v^*\in T_xM$ dual to
$v$.

The chirality operator \refe{Gam} coincides in this case (cf.
\cite[{\S}3.6]{BeGeVe}) with the Hodge $*$-operator. Hence, the usual grading
$\Lam{T^*M}=\Lam^{\text{even}}{T^*M}\oplus\Lam^{\text{odd}}{T^*M}$ is not the
natural grading in the sense of \refss{chir}. {\em We will always consider
$\Lam{T^*M}$ with the natural grading}. The positive and negative elements of
$\Gam(M,\Lam{T^*M})$, with respect to this grading, are called {\em self-dual
\/} and {\em anti-self-dual \/} differential forms respectively.

The action \refe{cv} is self-adjoint with respect to the metric on $\Lam{T^*M}$
defined by the Riemannian metric on $M$. The connection induced on $\Lam{T^*M}$
by the Levi-Civita connection on $T^*M$ is a Clifford connection. The Dirac
operator associated with this connection is equal to $d+d^*$ and is called the
{\em signature operator}, \cite[{\S}3.6]{BeGeVe}. If the dimension of $M$ is
divisible by four, then its index is equal to the signature of the manifold $M$.

\ssec{almcomp}{Almost complex manifolds}
Assume that $M$ is an almost complex manifold with an almost complex structure
$J:TM\to TM$. Then $J$ defines a structure of a complex vector bundle on the
tangent bundle $TM$. Let $h^{TM}$ be a Hermitian metric on $TM\otimes\CC$. The
real part $g^{TM}=\RE h^{TM}$ of $h^{TM}$ is a Riemannian metric on $M$. Note
also that $J$ defines an orientation on $M$.


Let $\Lam^q =\Lam^q(T^{0,1}M)^*$ denote the bundle of $(0,q)$-forms on $M$ and
set
\[
        \Lam^+ \ = \ \bigoplus_{q \text{ even}} \Lam^q, \quad
        \Lam^-  \ = \ \bigoplus_{q \text{ odd}} \Lam^q.
\]

Let $\lam^{1/2}$ be the square root of the complex line bundle $\lam=\det
T\ha{M}$ and let $\calS$ be the spinor bundle over $M$ associated to the
Riemannian metric $g^{TM}$. Although $\lam^{1/2}$ and $\calS$ are defined only
locally, unless $M$ is a spin manifold, it is well known (cf.
\cite[Appendix~D]{LawMic89}) that the products $\calS^\pm\otimes\lam^{1/2}$ are
globally defined and
\[
        \Lam^\pm \ = \ \calS^\pm\otimes\lam^{1/2}.
\]
It follows that $\Lam$ is a $\spin^c$ vector bundle over $M$ (cf. \refss{spin}).
In particular, the grading $\Lam=\Lam^+\oplus\Lam^-$ is natural.

The Clifford action of $C(M)$ on $\Lam$ may be described as follows: if
$f\in\Gam(M,T^*M)$ decomposes as $f=f\ha+f\ah$ with $f\ha\in\Gam(M,(T\ha M)^*)$
and $f\ah\in\Gam(M,(T\ah M)^*)$, then the Clifford action of $f$ on
$\alp\in\Gam(M,\Lam)$ equals
\begin{equation}\label{E:claction}
        c(f)\alp \ = \ \sqrt{2}\, \left( f\ah\wedge\alp  \ - \ \iot(f\ha)\, \alp\right).
\end{equation}
Here $\iot(f\ha)$ denotes the interior multiplication by the vector field
$(f\ha)^*\in T\ah M$ dual to the 1-form $f\ha$. This action is self-adjoint with
respect to the Hermitian structure on $\Lam$ defined by the Riemannian metric
$g^{TM}$ on $M$.

The Levi-Civita connection $\n^{TM}$ of $g^{TM}$ induces a Hermitian connection
on $\lam^{1/2}$ and on $\calS$. Let $\n^M$ be the product connection (cf.
\refss{clmodule}),
\[
        \n^M \ = \ \n^\calS\otimes 1 \ + \ 1\otimes\n^{\lam^{1/2}}.
\]
Then $\n^M$ is a well-defined Hermitian Clifford connection on the $\spin^c$
bundle $\Lam$. Hence, it gives rise to a self-adjoint $\spin^c$ Dirac operator.

More generally, assume that $\W$ is a Hermitian vector bundle over $M$ and let
$\n^\W$ be a Hermitian connection on $\W$. Consider the twisted Clifford module
$\E=\Lam\otimes\W$. The product connection $\n^\E=\n^{\Lam\otimes\W}$ determines
a Dirac operator $D:\Gam(M,\E)\to\Gam(M,\E)$.

\ssec{comp}{\ka manifolds}
If $(M,J,g^{TM})$ is a \ka manifold, then (cf. \cite[Proposition~3.67]{BeGeVe})
the Dirac operator defined in \refe{dir1} coincides with the Dolbeault-Dirac
operator
\begin{equation}\label{E:dirka}
        D \ = \ \sqrt2\, \left(\p+\p^*\right).
\end{equation}
Here $\p^*$ denotes the operator adjoint to $\p$ with respect to the
$L_2$-scalar product on $\A^{0,*}(M,\W)$. Hence, the restriction of the kernel
of $D$ to $\A^{0,i}(M,\W)$ is isomorphic to the cohomology $H^i(M,\O(\W))$ of
$M$ with coefficients in the sheaf of holomorphic section of $\W$.

\sec{main}{Vanishing theorems and their applications}

In this section we state the main theorems of the paper. The section is organized as
follows:

In \refss{lbundle}, we formulate our main result -- the vanishing theorem for the
half-kernel of a Dirac operator (\reft{main}).

In \refss{sofproof}, we briefly indicate the idea of the proof of \reft{main}.

In \refss{sign}, we apply this theorem to calculate the sign of the signature of a
vector bundle twisted by a high power of a line bundle.

In \refss{cvanish}, we refine \reft{main} for the case of a complex manifold. In
particular, we recover the Andreotti-Grauert vanishing theorem for a line bundle
with curvature of a mixed sign, cf. \cite{AndGr62,DemPetSch93}.

{}Finally, in \refss{acvanish}, we present an analogue of the Andreotti-Grauert
theorem for almost complex manifolds. This generalizes a result of Borthwick and
Uribe \cite{BoUr96}.

\ssec{lbundle}{Twisting by a line bundle. The vanishing theorem}
Suppose $\E$ is a self-adjoint Clifford module over $M$. Recall from
\refss{chir} that $\E=\E^+\oplus\E^-$ denotes the {\em natural \/} grading on
$\E$. Let $\L$ be a Hermitian line bundle over $M$, and let $\n^\L$ be a
Hermitian connection on $\L$. The connection $\nk$ (cf. \refe{nEW}) is a
Hermitian Clifford connection on the twisting Clifford module $\Ek$. Consider
the self-adjoint Dirac operator
\[
        D_k:\, \Gam(M,\E\otimes\Lk) \ \to \ \Gam(M,\E\otimes\Lk)
\]
associated to this connection and let $D_k^\pm$ denote the restriction of $D_k$
to the spaces $\Gam(M,\E^\pm\otimes\Lk)$.

The curvature $F^\L=(\n^\L)^2$ of the connection $\n^\L$ is an imaginary valued
closed 2-form on $M$. If it is non-degenerate, then $iF^\L$ is a symplectic form
on $M$ and, hence, defines an orientation of $M$.
Our main result is the following
\th{main}
Let $\E$ be a self-adjoint Clifford module over a compact oriented
even-dimensional Riemannian manifold $M$. Let $\n^\E,\L$, $\n^\L,D_k$ be as
above. Assume that the curvature $F^\L=(\n^\L)^2$ of the connection $\n^\L$ is
non-degenerate at all points of $M$. If the orientation defined by the
symplectic form $iF^\L$ coincides with the original orientation of $M$, then
\begin{equation}\label{E:main}
        \Ker D^-_k \ = \ 0 \qquad \mbox{for} \qquad k\gg 0.
\end{equation}
Otherwise, $\Ker D^+_k = 0$ for $k\gg 0$.
\eth
The theorem is a generalization of a vanishing theorem of Borthwick and Uribe
\cite{BoUr96}, who considered the case where $M$ is an almost \ka manifold, $D$
 is a  $\spin^c$-Dirac operator and  $\L$  is a positive line bundle
over $M$.

The theorem is proven in \refs{prmain}. Here we only explain the main ideas of the
proof.
%
\ssec{sofproof}{The scheme of the proof}
Our proof of Theorem~\ref{T:main} follows the lines of \cite{BoUr96}. It is
based on an estimate from below on the large $k$ behavior of the square $D_k^2$
 of the Dirac operator. Using this estimate we show that, if the orientation
defined by $iF^\L$ coincides with (resp. is opposite to) the given orientation
of $M$, then, for large $k$, the restriction of $D_k^2$ to $\E^-\otimes\Lk$
(resp. to $\E^+\otimes\Lk$) is a strictly positive operator and, hence, has no
kernel.

This estimate for $D_k^2$ is obtained in two steps. First we use the
Lichnerowicz formula (cf. \refss{lichn}) to compare $D_k^2$ with the {\em metric
Laplacian \/} $\Del_k=\n^{\E\otimes\Lk}(\n^{\E\otimes\Lk})^*$.

Then it remains to study the large $k$ behavior of the metric Laplacian
$\Del_k$. This is done in \refs{lapl}. In fact, the estimate which we need is
essentially obtained in \cite{BoUr96,GuiUr88}. Roughly speaking it says that
$\Del_k$ grows linearly in $k$.

The proof of the estimate for $\Del_k$ also consists of two steps. {}First we
consider the principal bundle $\Z\to M$ associated to the vector bundle
$\E\otimes\L$, and construct a differential operator $\tilDel$ ({\em horizontal
Laplacian}) on $\Z$, such that the operator $\Del_k$ is ``equivalent" to a
restriction of $\tilDel$ on a certain subspace of the space of $L_2$-functions
on the total space of $\Z$. Then we apply the {\em a priory \/} Melin estimates
\cite{Melin71} (see also \cite[Theorem~22.3.3]{Horm3}) to the operator
$\tilDel$.
\rem{exest}
It would be very interesting to obtain an effective estimates of a
minimal value of $k$ which satisfies \refe{main} at least for the
simplest cases (say, when $\E$ is a spinor bundle over a spin manifold
$M$). Unfortunately, such an estimate can not be obtained using our
method. This is because the Melin inequalities \cite{Melin71},
\cite[Theorem~22.3.3]{Horm3} (see also \refss{melin}), used in our
proof, contain a constant $C$, which can not be estimated effectively.
\erem

We will now discuss applications and refinements of \reft{main}. In particular,
we will see that \reft{main} may be considered as a generalization of the
vanishing theorems of Kodaira, Andreotti-Grauert \cite{AndGr62} and
Borthwick-Uribe \cite{BoUr96}.

\ssec{sign}{The signature operator}
Recall from \refss{difforms} that, for any oriented even-\linebreak{dimensional}
Riemannian manifold $M$, the exterior algebra $\Lam{T^*M}$ of the cotangent
bundle is a self-adjoint Clifford module. The connection induced on $\Lam{T^*M}$
by a Levi-Civita connection on $T^*M$ is a Hermitian Clifford connection and the
Dirac operator associated to this connection is the signature operator $d+d^*$.

Consider a twisted Clifford module $\E=\Lam{T^*M}\otimes\W$, where $\W$ is a
Hermitian vector bundle over $M$ endowed with a Hermitian connection $\n^\W$.

Let $\L$ be a Hermitian line bundle over $M$ and let $\n^\L$ be an Hermitian
connection on $\L$. The space $\Gam(M,\E\otimes\Lk)$ of sections of the twisted
Clifford module $\E\otimes\Lk$ coincides with the space $\A^*(M,\W\otimes\Lk)$
 of differential forms on  $M$  with values in  $\W\otimes\Lk$. The
positive and negative elements of $\A^*(M,\W\otimes\Lk)$ with respect to the
natural grading are called the {\em self-dual \/} and the {\em anti-self-dual
\/} differential forms respectively.

Let $D_k:\A^*(M,\W\otimes\Lk)\to\A^*(M,\W\otimes\Lk)$ denote the Dirac operator
corresponding to the tensor product connection $\n^{\W\otimes\Lk}$ on
$\W\otimes\Lk$. Then
\begin{equation}\label{E:derham}
        D_k \ = \ \n^{\W\otimes\Lk} \ + \ \left(\n^{\W\otimes\Lk}\right)^*,
\end{equation}
where $\left(\n^{\W\otimes\Lk}\right)^*$ denotes the adjoint of
$\n^{\W\otimes\Lk}$ with respect to the $L_2$-scalar product on $\W\otimes\Lk$.
The operator \refe{derham} is called the {\em signature operator \/} of the
bundle $\W\otimes\Lk$.

Let $D_k^+$ and $D_k^-$ denote the restrictions of $D_k$ on the spaces of
self-dual and anti-self-dual differential forms respectively.

As an immediate consequence of \reft{main}, we obtain the following
\th{derham}
Suppose $M$ is a compact oriented even-dimensional Riemannian manifold. Let $\W,
\n^\W, \L, \n^\L, D_k^\pm$ be as above. Assume that the curvature $F^\L=(\n^\L)^2$
of the connection $\n^\L$ is non-degenerate at any point $x\in M$. If the
orientation defined by $iF^\L$ coincides with the given orientation of $M$, then
$\Ker D^-_k=0$, for $k\gg0$.

Otherwise, $\Ker D^+_k=0$, for $k\gg0$.
\eth

The index
\[
        \ind D_k \ = \ \dim\Ker D_k^+ \ - \ \dim\Ker D_k^-
\]
of the Dirac operator $D_k$ is called the {\em signature \/} of the bundle
$\W\otimes\Lk$ and is denoted by $\sign(\W\otimes\Lk)$. It depends only on the
manifold $M$, its orientation and the bundle $\W\otimes\Lk$ (but not on the
choice of Riemannian metric on $M$ and of Hermitian structures and connections
on the bundles $\W,\L$). If the bundles $\W$ and $\L$ are trivial, then it
coincides with the usual signature of the manifold $M$.

From \reft{derham}, we obtain the following
\cor{derham}
Let $\W$ and $\L$ be respectively a vector and a line bundles over a compact
oriented even-dimensional Riemannian manifold $M$. Suppose that, for some
Hermitian metric on $\L$, there exist a Hermitian connection, whose curvature
$F^\L$ is non-degenerate at any point of $M$. If the orientation defined by the
symplectic form $iF^\L$ coincides with the given orientation of $M$, then
\[
        \sign(\W\otimes\Lk) \ \ge \ 0 \qquad \mbox{for} \qquad k\gg 0.
\]
Otherwise, $\sign(\W\otimes\Lk)\le0$ for $k\gg0$.
\ecor

\ssec{cvanish}{Complex manifolds. The Andreotti-Grauert theorem}
Suppose $M$ is a compact complex manifold, $\W$ is a holomorphic vector bundle
over $M$ and $\L$ is a holomorphic line bundles over $M$. Fix a Hermitian metric
$h^\L$ on $\L$ and let $\n^\L$ be the {\em Chern connection \/} on $\L$, i.e.,
the unique holomorphic connection which preserves the Hermitian metric. The
curvature $F^\L$ of $\n^\L$ is a $(1,1)$-form which is called the {\em curvature
form of the Hermitian metric $h^\L$}.

The orientation condition of \reft{main} may be reformulated as follows. Let
$(z^1\nek z^n)$ be complex coordinates in the neighborhood of a point $x\in M$.
The curvature $F^\L$ may be written as
\[
        iF^\L \ = \ \frac{i}2\sum_{i,j} F_{ij}dz^i\wedge d\oz^j.
\]
Denote by $q$ the number of negative eigenvalues of the matrix $\{F_{ij}\}$.
Clearly, the number $q$ is independent of the choice of the coordinates. We will
refer to this number as the {\em number of negative eigenvalues of the curvature
$F^\L$ at the point $x$}. Then the orientation defined by the symplectic form
$iF^\L$ coincides with the complex orientation of $M$ if and only if $q$ is
even.

A small variation of the method used in the proof of \reft{main} allows to get a
more precise result which depends not only on the parity of $q$ but on $q$
itself. In this way we obtain a new proof of the following vanishing theorem of
Andreotti and Grauert \cite{AndGr62,DemPetSch93}
%
\begin{Thm}[{\textbf{Andreotti-Grauert}}]\label{T:AG}
Let $M$ be a compact complex manifold and let $\L$ be a holomorphic line bundle
over $M$. Assume that $\L$ carries a Hermitian metric whose curvature form
$F^\L$ has at least $q$ negative and at least $p$ positive eigenvalues at any
point $x\in M$. Then, for any holomorphic vector bundle $\W$ over $M$, the
cohomology $H^j(M,\O(\W\otimes\Lk))$ with coefficients in the sheaf of
holomorphic sections of $\W\otimes\Lk$ vanish for $j\not=q, q+1\nek n-p$ and
$k\gg0$.
\end{Thm}
The proof is given in \refss{prAG}.
%
In contrary to \reft{main}, the curvature $F^\L$ in \reft{AG} needs not be
non-degenerate. If $F^\L$ is non-degenerate, then the number $q$ of negative
eigenvalues of $F^\L$ does not depend on the point $x\in M$. Then we obtain the
following
\cor{AG}
If, in the conditions of \reft{AG}, the curvature $F^\L$ is non-degenerate and
has exactly $q$ negative eigenvalues at any point $x\in M$, then
$H^j(M,\O(\W\otimes\Lk))$ vanishes for any $j\not=q$ and $k\gg0$.
\ecor
Note that, if $\L$ is a positive line bundle, \refc{AG} reduces to the classical
Kodaira vanishing theorem (cf., for example, \cite[Theorem~3.72(2)]{BeGeVe}).
\rem{2}
a. \ It is interesting to compare \refc{AG} with \reft{main} for the case when
$M$ is a \ka manifold. In this case the Dirac operator $D_k$ is equal to the
Dolbeault-Dirac operator \refe{dirka}. Hence (cf. \refss{comp}), \reft{main}
implies that $H^j(M,\O(\W\otimes\Lk))$ vanishes when the parity of $j$ is not
equal to the parity of $q$. \refc{AG} refines this result.

b. \ If $M$ is not a \ka manifold, then the Dirac operator $D_k$ defined by
\refe{dir1} is not equal to the Dolbeault-Dirac operator, and the kernel of
$D_k$ is not isomorphic to the cohomology $H^*(M,\O(\W\otimes\Lk))$. However, we
show in \refs{prAG} that the operators $D_k$ and $\sqrt(\p+\p^*)$ have the same
asymptotic as $k\to\infty$. Then the vanishing of the kernel of $D_k$ implies
the vanishing of the cohomology $H^j(M,\O(\W\otimes\Lk))$.

c. \ In \reft{AG} the bundle $\W$ can be replaced by an arbitrary coherent sheaf
 $\calF$. This follows from \reft{AG} by a standard technique using a
resolution of $\calF$ by locally free sheaves (see, for example,
\cite[Ch.~5]{ShSom85} for similar arguments).
\erem

\ssec{acvanish}{Andreotti-Grauert-type theorem for almost complex manifolds}
In this section we refine \reft{main} assuming that $M$ is endowed with an
almost complex structure $J$ such that the curvature $F^\L$ of $\L$ is a $(1,1)$
 form on  $M$  with respect to  $J$. In other words, we assume that, for
any $x\in M$ and any basis $(e^1\nek e^n)$ of the holomorphic cotangent space
$(T\ha{M})^*$, one has
\[
        iF^\L \ = \ \frac{i}2\sum_{i,j} F_{ij}e^i\wedge \oe^j.
\]

This section generalizes a result of Borthwick and Uribe \cite{BoUr96}.

We denote by $q$ the number of negative eigenvalues of the matrix $\{F_{ij}\}$.
As in \refss{cvanish}, the orientation of $M$ defined by the symplectic form
$iF^\L$ depends only on the parity of $q$. It coincides with the orientation
defined by $J$ if and only if $q$ is even.

We will use the notation of \refss{almcomp}. In particular, $\Lam=\Lam(T\ah
M)^*$ denotes the bundle of $(0,*)$-forms on $M$ and $\W$ is a Hermitian vector
bundle over $M$. Then $\E=\Lam\otimes\W$ is a self-adjoint Clifford module
endowed with a Hermitian Clifford connection $\n^\E$. The space
$\Gam(M,\E\otimes\Lk)$ of sections of the twisted Clifford module $\E\otimes\Lk$
coincides with the space $\A^{0,*}(M,\W\otimes\Lk)$ of differential forms of
type $(0,*)$ with values in $\W\otimes\Lk$. Let
\[
        D_k:\, \A^{0,*}(M,\W\otimes\Lk) \ \to \ \A^{0,*}(M,\W\otimes\Lk)
\]
denote the Dirac operator corresponding to the tensor product connection on
$\W\otimes\Lk$.

{}For a form $\alp\in \A^{0,*}(M,\W\otimes\Lk)$, we denote by $\|\alp\|$ its
$L_2$-norm and by $\alp_i$ its component in $\A^{0,i}(M,\W\otimes\Lk)$.
%
\th{UB+}
Assume that the matrix $\{F_{ij}\}$ has at least $q$ negative and at least $p$
positive eigenvalues at any point $x\in M$. Then there exists a sequence $\eps_1,\eps_2,\ldots$ convergent to zero, such that for
any $k\gg 0$ and any $\alp\in\Ker D_k^2$ one has
\[
        \|\alp_j\| \ \le \  \eps_k\|\alp\|, \qquad \mbox{for} \quad
                                                        j\not=q,q+1\nek n-p.
\]

In particular, if the form $F^\L$ is non-degenerate and $q$ is the number of
negative eigenvalues of $\{F_{ij}\}$ (which is independent of $x\in M$), then
there exists a sequence $\tileps_{1},\tileps_{2},\ldots$, convergent to zero,
such that $\alp\in\Ker D_k$ implies
\[
         \|\alp \ - \ \alp_q\| \ \le \ \tileps_{k}\|\alp_q\|.
\]
\eth
\reft{UB+} is proven in \refss{prUB}. The main ingredient of the proof is the
following estimate on $D_k$, which also has an independent interest:
%
\prop{estDk}
If the matrix $\{F_{ij}\}$ has at least $q$ negative and at least $p$ positive
eigenvalues at any point $x\in M$, then there exists a constant $C>0$, such
that
\[
        \|D_k\, \alp\| \ \ge \ Ck^{1/2}\, \|\alp\|,
\]
for any $k\gg0, \ j\not=q,q+1\nek n-p$ and $\alp\in
\A^{0,j}(M,\W\otimes\Lk)$.
\eprop
The proof is given in \refss{prestDk}.
\rem{UB}
a. \ For the case when $\L$ is a positive line bundle, the Riemannian metric on
$M$ is almost \ka and $\W$ is a trivial line bundle, \reft{UB+} was established
by Borthwick and Uribe \cite[Theorem~2.3]{BoUr96}.

b. \ \reft{UB+} implies that, if $F^\L$ is non-degenerate, then $\Ker D_k$ is
dominated by the component of degree $q$. If $\alp\in \Gam(M,\E^-)$ (resp.
$\alp\in\Gam(M,\E^+)$) and $q$ is even (resp. odd) then $\alp_q=0$. So, we
obtain the vanishing result of \reft{main} for the case when $M$ is almost
complex and $F^\L$ is a $(1,1)$ form.

c. \ \reft{UB+} is an analogue of \reft{AG}. Of course, the cohomology
$H^j(M,\O(\W\otimes\Lk))$ is not defined if $J$ is not integrable. Moreover, the
square $D^2_k$ of the Dirac operator does not preserve the $\ZZ$-grading on
$\A^{0,*}(M,\W\otimes\Lk)$. Hence, one can not hope that the kernel of $D_k$
belongs to $\oplus_{j=q}^{n-p}\A^{0,j}(M,\W\otimes\Lk)$. However, \reft{UB+}
shows, that for any $\alp\in\Ker D_k$, ``most of the norm" of $\alp$ is
concentrated in $\oplus_{j=q}^{n-p}\A^{0,j}(M,\W\otimes\Lk)$.
\erem


\sec{prmain}{Proof of the vanishing theorem for the half-kernel of a Dirac operator}

In this section we present a proof of \reft{main} based on
Propositions~\ref{P:D-D} and \ref{P:lapl}, which will be proved in the
following sections.

The idea of the proof is to study the large $k$ behavior of the square $D_k^2$
of the Dirac operator.

\ssec{tJ}{The operator $\tJ$}
We need some additional definitions. Recall that $F^\L$ denotes the curvature of
the connection $\n^\L$. In this subsection we do not assume that $F^\L$ is
non-degenerate. For $x\in M$, define the skew-symmetric linear map
$\tJ_x:T_xM\to T_xM$ by the formula
\[
        iF^\L(v,w) \ = \ g^{TM}(v,\tJ_xw), \qquad v,w\in T_xM.
\]
The eigenvalues of $\tJ_x$ are purely imaginary. Note that, in general, $\tJ$ is
not an almost complex structure on $M$.

Define
\begin{equation}\label{E:tau}
        \tau(x) \ = \ \Tr^+ \tJ_x \ := \ \mu_1+\cdots+\mu_l, \qquad
        m(x)    \ = \ \min_{j}\, \mu_j(x).
\end{equation}
where $i\mu_j, \ j=1\nek l$ are the eigenvalues of $\tJ_x$ for which $\mu_j>0$.
Note that $m(x)=0$ if and only if the curvature $F^\L$ vanishes at the point
$x\in M$.

\ssec{lapl}{Estimate on $D_k^2$}
Our estimate on the square $D_k^2$ of the Dirac operator is obtained in two
steps: first we compare it to the {\em metric Laplacian}
\[
        \Del_k \ := \ (\n^{\E\otimes\Lk})^*\, \n^{\E\otimes\Lk},
\]
and then we estimate the large $k$ behavior of $\Del_k$. These two steps are the
subject of the following two propositions.
%
\prop{D-D}
Supposed that the differential form $F^\L$ is non-degenerate. If the orientation
defined on $M$ by the symplectic form $iF^\L$ coincides with (resp. is opposite
to) the given orientation of $M$, then there exists a constant $C$ such that,
for any $s\in\Gam(M,\E^-\otimes\Lk)$ (resp. for any $s\in
\Gam(M,\E^+\otimes\Lk)$), one has an estimate
\[
        \big\langle\, (D^2_k - \Del_k)\, s,\, s\, \big\rangle \ \ge \
                -k\,\big\langle\, (\tau(x)-2m(x))\, s,\, s\, \big\rangle - C\, \|s\|^2.
\]
Here $\<\cdot,\cdot\>$ denotes the $L_2$-scalar product on the space of sections
and $\|\cdot\|$ is the norm corresponding to this scalar product.
\eprop
The proposition is proven in \refss{lichn} using the Lichnerowicz formula \refe{lichn}.

In the next proposition we do not assume that $F^\L$ is non-degenerate.
\prop{lapl}
Suppose that $F^\L$ does not vanish at any point $x\in M$. For any $\eps>0$,
there exists a constant $C_\eps$ such that, for any $k\in\ZZ$ and any section
$s$ of the bundle $\E\otimes\Lk$,
\eq{lapl}
        \<\Del_k\, s,s\> \ \ge \ k\, \<(\tau(x)-\eps) s,s\> \ - \ C_\eps\, \|s\|^2.
\end{equation}
\eprop
\refp{lapl} is essentially proven in \cite[Theorem~2.1]{BoUr96}. The only
difference is that we do not assume that the curvature $F^\L$ has a constant
rank. This forces us to use the original Melin inequality \cite{Melin71} (see
also \cite[Theorem~22.3.3]{Horm3}) and not the H\"ormander refinement of this
inequality \cite[Theorem~22.3.2]{Horm3}. That is the reason that $\eps\not=0$
appears in \refe{lapl}. Note, \cite{BoUr96}, that if $F^\L$ has constant rank,
then \refp{lapl} is valid for $\eps=0$.

\refp{lapl} is proven in \refs{lapl}.

\ssec{prmain}{Proof of \reft{main}}
Assume that the orientation defined by $iF^\L$ coincides with the given
orientation of $M$ and $s\in\Gam(M,\E^-\otimes\L)$, or that the orientation
defined by $iF^\L$ is opposite to the given orientation of $M$ and
$s\in\Gam(M,\E^+\otimes\L)$. By \refp{D-D},
\eq{D>D}
        \<D_k^2\, s,s\> \ \ge \ \< \Del_k\, s,s\> \ - \
        k\,\big\langle\, (\tau(x)-2m(x))\, s,\, s\big\rangle \ - \ C\, \|s\|^2.
\end{equation}
Choose
\[
        0 \ < \ \eps \ < \ 2\, \min_{x\in M} m(x)
\]
and set
\[
        C'=2\min_{x\in M} m(x)-\eps>0.
\]
Since the metric Laplacian $\Del_k$ is a non-negative operator, it follows from
\refe{lapl} and \refe{D>D} that
\[
         \<D_k^2\, s,\, s\> \ \ge \ kC'\|s\|^2 \ - \ (C+C_\eps)\, \|s\|^2.
\]
Thus, for $k> (C+C_\eps)/C'$, we have $\<D_k^2\, s,s\>>0$. Hence, $D_k s\not=0$.
\hfill$\square$

\sec{prUB}{Proof of the vanishing theorem for almost complex manifolds}

In this section we prove \reft{UB+} and \refp{estDk}. The proof is very similar
to the proof of \reft{main} (cf. \refs{prmain}). It is based on \refp{lapl} and
the following refinement of \refp{D-D}:

\prop{D-Dac}
Assume that the matrix $\{F_{ij}\}$ (cf. \refss{acvanish}) has at least $q$
negative eigenvalues at any point $x\in M$. For any $x\in M$, we denote by
$m_q(x)>0$ the minimal positive number, such that at least $q$ of the
eigenvalues of $\{F_{ij}\}$ do not exceed $-m_q$. Then there exists a constant
 $C$  such that
\[
        \big\langle\, (D^2_k - \Del_k)\, \alp,\, \, \alp\big\rangle \ \ge \
                -k\, \big\langle\, (\tau(x)-2m_q(x))\, \alp,\, \alp\, \big\rangle
                \  - \ C\, \|\alp\|^2
\]
for any $j=0\nek q-1$ and any $\alp\in \A^{0,j}(M,\W\otimes\Lk)$.
\eprop
The proposition is proven in \refss{prD-Dac}.

\ssec{prestDk}{Proof of \refp{estDk}}
Choose $0<\eps<2\min_{x\in M}\, m_q(x)$ and set
\[
        C' \ = \ 2\min_{x\in M}\, m_q(x)-\eps.
\]
{}Fix $j=0\nek q-1$ and let $\alp\in \A^{0,j}(M,\W\otimes\Lk)$. Since the metric
Laplacian $\Del_k$ is a non-negative operator, it follows from
Propositions~\ref{P:lapl} and \ref{P:D-Dac}, that
\[
        \< D_k^2\, \alp,\alp\> \ \ge \ kC'\, \|\alp\|^2 \ - \ (C+C_\eps)\, \|\alp\|^2.
\]
Hence, for any $k>2(C+C_\eps)/C'$, we have
\[
        \|D_k\, \alp\|^2 \ = \ \<D_k^2\, \alp,\, \alp\, \> \ \ge \
                \frac{k C'}2\, \|\alp\|^2.
\]
This proves \refp{estDk} for $j=0\nek q-1$. The statement for $j=n-p+1\nek n$
may be proven by a verbatim repetition of the above arguments, using a natural
analogue of \refp{D-Dac}. (Alternatively, the statement for $j=n-p+1\nek n$ may
be obtained as a formal consequence of the statement for $j=0\nek q-1$ by
considering $M$ with an opposite almost complex structure).
\hfill$\square$

\ssec{gr}{}
If the manifold $M$ is not K\"ahler, then the operator $D_k^2$ does not preserve
the $\ZZ$-grading on $\A^{0,*}(M,\W\otimes\Lk)$. However, the next proposition
shows that the {\em mixed degree operator} $\alp_i\mapsto(D_k^2\alp_i)_j$ is, in
a certain sense, small.
\prop{gr}
There exists a sequence $\del_{1}, \del_{2}\nek$ such that 
$\lim_{k\to\infty}\del_{k}=0$ and
\[
        |\<D_k^2\, \alp, \bet\>| \ \le \del_{k}\, \<D_k^2\, \alp, \alp\> \ + \
          \del_{k}\, \<D_k^2\, \bet, \bet\> \ + \
             \del_{k}k\, \|\alp\|^{2} \ + \  \del_{k}k\, \|\bet\|^{2},
\]
for any $i\not=j$ and any $\alp\in \A^{0,i}(M,\W\otimes\Lk), 
\bet\in \A^{0,j}(M,\W\otimes\Lk)$.
\eprop
The proof of the proposition, based on the Lichnerowicz formula, is given in
\refss{prgr}.

\ssec{prUB}{Proof of \reft{UB+}}
Let $\alp\in \Ker D_k$ and fix $j\not\in q,q+1\nek n-p$. Set $\bet=\alp-\alp_j$.
Then
\[
        0 \ = \ \|D_k\, \alp\|^2 \ = \ \|D_k\, \alp_j\|^2  \ + \
               2\RE\, \<\, D_k\, \alp_j,\, D_k\, \bet\, \> \ + \
                  \|D_k\, \bet\|^2.
\]
Hence, it follows from \refp{gr}, that
\begin{equation}\label{E:alpbet}
        0 \ \ge \ (1-2\del_k)\, \|D_k\, \alp_j\|^2 \ + \
           (1-2\del_k)\, \|D_k\, \bet\|^2 \ - \
                  2\del_{k}k\, \|\alp_i\|^{2} \ - \ 2\del_{k}k\, \|\bet\|^{2}.
\end{equation}
If we assume now that $k$ is large enough, so that $1-2\del_k>0$, then we obtain
from \refe{alpbet} and \refp{estDk} that
\[
        \big((1-2\del_k)C^2k -2\del_kk\big)\,  \|\alp_i\|^{2} \ \le \
                        2\del_{k}k\, \|\bet\|^{2} \ \le \ 2\del_{k}k\, \|\alp\|^{2}.
\]
Thus
\[
        \|\alp_i\|^{2} \ \le \
                \frac{2\del_{k}}{(1-2\del_k)C^2 -2\del_k}\, \|\alp\|^{2}.
\]
Hence, \reft{UB+} holds with
$\eps_k=\sqrt{\frac{2\del_{k}}{(1-2\del_k)C^2-2\del_k}}$.
\hfill$\square$

\sec{prAG}{Proof of the Andreotti-Grauert theorem}

In this section we use the results of \refss{acvanish} in order to get a new proof of
the Andreotti-Grauert theorem (\reft{AG}).

Note first, that, if the manifold $M$ is K\"ahler, then the Andreotti-Grauert
theorem follows directly from \reft{UB+}. Indeed, in this case the Dirac
operator $D_k$ is equal to the Dolbeault-Dirac operator $\sqrt2(\p+\p^*)$.
Hence, the restriction of the kernel of $D_k$ to $\A^{0,j}(M,\O(\W\otimes\Lk))$
 is isomorphic to the cohomology  $H^j(M,\O(\W\otimes\Lk))$.

In general, $D_k\not=\sqrt2(\p+\p^*)$. However, the following proposition shows
that those two operators have the same ``large $k$ behavior".

Recall from \refss{acvanish} the notation
\[
        \E \ = \ \Lam(T^{0,1}M)\otimes\W.
\]

\prop{D-dd}
There exists a bundle map $A\in\End(\E)\subset \End(\E\otimes\Lk)$, independent
of $k$, such that
\begin{equation}\label{E:D-dd}
        \sqrt2\, \left(\p+\p^*\right) \ = \ D_k \ + \ A.
\end{equation}
\eprop
\prf
Choose a holomorphic section $e(x)$ of $\L$ over an open set $U\subset\L$. It
defines a section $e^k(x)$ of $\Lk$ over $U$ and, hence, a holomorphic
trivialization
\begin{equation}\label{E:triv}
        U\times\CC \ \overset{\sim}{\longrightarrow} \
                \Lk, \qquad (x,\phi)\mapsto \phi\cdot e^k(x)\in \Lk
\end{equation}
of the bundle $\Lk$ over $U$. Similarly, the bundles $\W$ and $\W\otimes\Lk$ may
be identified over $U$ by the formula
\begin{equation}\label{E:WLk=W}
        w \ \mapsto \ w\otimes e^k.
\end{equation}

Let $h^\L$ and $h^\W$ denote the Hermitian fiberwise metrics on the bundles $\L$
 and  $\W$  respectively. Let  $h^{\W\otimes\Lk}$  denote the
Hermitian metric on $\W\otimes\Lk$ induced by the metrics $h^\L, h^\W$. Set
\[
        f(x)\ := \ |e(x)|^2, \qquad x\in U,
\]
where $|\cdot|$ denotes the norm defined by the metric $h^\L$. Under the
isomorphism \refe{WLk=W} the metric $h^{\W\otimes\Lk}$ corresponds to the metric
\begin{equation}\label{E:hk}
        h_k(\cdot,\cdot) \ = \ f^k\, h^\W(\cdot,\cdot)
\end{equation}
on $\W$.

By \cite[p.~137]{BeGeVe}, the connection $\n^\L$ on $\L$ corresponds under the
trivialization \refe{triv} to the operator
\[
        \Gam(U,\CC) \ \to \ \Gam(U,T^*U\otimes\CC);     \qquad
                        s\mapsto ds+kf^{-1}\d f\wedge s.
\]
Similarly, the connection on $\E\otimes\Lk= \Lam(T\ah M)^*\otimes\W\otimes\Lk$
corresponds under the isomorphism \refe{WLk=W} to the connection
\[
        \n_k:\, \alp \ \mapsto \ \n^\E\alp \ + \ kf^{-1}\d f\wedge\alp, \qquad
                \alp\in \Gam(U,\Lam(T\ah U)^*\otimes\W|_U)
\]
on $\E|_U$. It follows now from \refe{claction} and \refe{dir1} that the Dirac
operator $D_k$ corresponds under \refe{WLk=W} to the operator
\begin{equation}\label{E:tilDk}
        \tilD_k:\, \alp \ \mapsto \ D_0\alp \ - \ \sqrt2kf^{-1}\iot(\d f)\alp, \qquad
                        \alp\in \A^{0,*}(U,\W|_U).
\end{equation}
Here $\iot(\d f)$ denotes the contraction with the vector field $(\d f)^*\in
T\ah{M}$ dual to the 1-form $\d{f}$, and $D_0$ stands for the Dirac operator on
the bundle $\E=\E\otimes\L^0$.

Let $\p_k^*:\A^{0,*}(U,\W|_U) \to \A^{0,*-1}(U,\W|_U)$ denote the adjoint of the
operator $\p$ with respect to the scalar product on $\A^{0,*}(U,\W|_U)$
determined by the Hermitian metric $h_k$ on $\W$ and the Riemannian metric on
$M$. Then, it follows from \refe{hk}, that
\begin{equation}\label{E:pk}
        \p_k^* \ = \ \p \ + \ kf^{-1}\iot(\d f).
\end{equation}
%
By \refe{tilDk} and \refe{pk}, we obtain
\[
        \sqrt2\, \left(\p+\p^*_k\right) \ - \ \tilD_k \ =  \
                \sqrt2\, \left(\p+\p^*_0\right) \ - \ D_0.
\]
Set $A=\sqrt2(\p+\p^*_0)-D_0$. By \cite[Lemma~5.5]{Duis96}, $A$ is a zero order
operator, i.e. $A\in\End(\E)$ (note that our definition of the Clifford action
on $\Lam(T\ah{M})^*$ and, hence, of the Dirac operator defers from \cite{Duis96}
by a factor of $\sqrt2$).
\eprf

\ssec{prAG}{Proof of \reft{AG}}
Let $A\in\End(\E)$ be the operator defined in \refp{D-dd} and let
\[
        \|A\| \ = \ \sup_{\|\alp\|=1} \, \|A\alp\|, \qquad
                        \alp\in \A^{0,*}(M,\E\otimes\Lk)
\]
be the $L_2$-norm of the operator $A:\A^{0,*}(M,\E\otimes\Lk)\to
\A^{0,*}(M,\E\otimes\Lk)$. By \refp{estDk}, there exists a constant  $C>0$  such
that
\[
        \|D_k\, \alp\| \ \ge \ Ck^{1/2}\, \|\alp\|,
\]
for any $k\gg0, \ j\not=q,q+1\nek n-p$ and $\alp\in
\A^{0,j}(M,\W\otimes\Lk)$. Then, if  $k> \|A\|^2/C^2$, we have
\[
        \|\sqrt2(\p+\p^*)\alp\| \ = \ \|(D_k+A)\alp\| \ \ge \
        \|D_k\alp\| \ - \ \|A\|\, \|\alp\|
        \ge \ \left(Ck^{1/2}-\|A\|\right)\, \|\alp\| \ > \ 0,
\]
for any $j\not=q,q+1\nek n-p$ and $0\not=\alp\in
\A^{0,j}(M,\W\otimes\Lk)$. Hence, the restriction of the kernel of the
Dolbeault-Dirac operator to the space $\A^{0,j}(M,\E\otimes\Lk)$ vanishes for
$j\not=q,q+1\nek n-p$.
\hfill$\square$

\sec{lichn}{The Lichnerowicz formula. Proof of Propositions~\ref{P:D-D},
         \ref{P:D-Dac} and \ref{P:gr}}

In this section we use the Lichnerowicz formula (cf. \refss{lichn}) to prove the
Propositions~\ref{P:D-D}, \ref{P:D-Dac} and \ref{P:gr}.

Before formulating the Lichnerowicz formula, we need some more information about
the Clifford modules and Clifford connections (cf. \cite[Section~3.3]{BeGeVe}).

\ssec{quant}{The symbol map and the quantization map}
Recall from \refss{difforms} that the exterior algebra $\Lam{T^*M}$ has a
natural structure of a self-adjoint Clifford module.

The Clifford bundle $C(M)$ is isomorphic to $\Lam{T^*M}$ as a bundle of vector
spaces. The isomorphism is given by the {\em symbol map}
\[
        \sig: \,  C(M) \ \to \Lam T^*M, \qquad \sig: \, a \ \mapsto c(a)1.
\]
The inverse of $\sig$ is called the {\em quantization map \/} and is denoted by
$\bfc$. If $e^1\nek e^{2n}$ is an orthonormal basis of $T^*_xM$, then (cf.
\cite[Proposition~3.5]{BeGeVe})
\[
        \bfc(e^{i_1}\wedge\dots\wedge e^{i_k}) \ = \ c(e^{i_1})\cdots c(e^{i_k}).
\]
Note that $\sig$ {\em is not \/} a map of algebras, i.e., $\sig(ab)\not=
\sig(a)\sig(b)$.

Assume now that $\E$ is a Clifford module. The composition of the quantization
map with the Clifford action of $C(M)$ on $\E$ defines a map $\bfc:\Lam{T^*M}\to
\End(\E)$. Though this map does not define the action of the exterior algebra on
$\E$ (i.e. $\bfc(ab)\not=\bfc(a)\bfc(b)$) it plays an important role in
differential geometry.

Let $\A(M)=\Gam(M,\Lam{T^*M})$ denote the space of smooth sections of
$\Lam{T^*M}$, i.e., the space of smooth differential forms on $M$. The
quantization map induces an isomorphism between $\A(M)$ and the space of
sections of $C(M)$. More generally, for any bundle $\E$ over $M$, there is an
isomorphism
\begin{equation}\label{E:AE=GCE}
        \A(M,\E) \ \cong \ \Gam(M,C(M)\otimes \E)
\end{equation}
between the space of differential forms on $M$ with values in $\E$ and the space
of smooth sections of the tensor product $C(M)\otimes \E$. Combining this
isomorphism with the Clifford action
\[
       c:\,  \Gam(M,C(M)\otimes \E) \ {\to} \ \Gam(M,\E),
\]
we obtain a map $\bfc:\, \A(M,\E) \to \E$. Similarly, we have a map
\begin{equation}\label{E:bfc}
        \bfc:\, \A(M,\End(\E)) \ \to \ \End(\E).
\end{equation}
In this paper we will be especially interested in the restriction of the above map to
the space of 2-forms. In this case the following formula is useful
\begin{equation}\label{E:bfc2}
        \bfc(F) \ = \ \sum_{i<j}\, F(e_i,e_j)\, c(e^i)\, c(e^j), \qquad
                        F\in \A^2(M,\End(\E)),
\end{equation}
where $(e_1\nek e_{2n})$ is an orthonormal frame of the tangent space to $M$,
and $(e^1\nek e^{2n})$ is the dual frame of the cotangent space.

\ssec{curv}{The curvature of a Clifford connection}
Let $\n^\E$ be a Clifford connection on a Clifford module $\E$ and let
\[
        F^\E \ = \ (\n^\E)^2 \ \in \ \A^2(M,\End(\E))
\]
denote the curvature of $\n^\E$.

Let $\EndC(\E)$ denote the bundle of endomorphisms of $\E$ commuting with the
action of the Clifford bundle $C(M)$. Then the bundle $\End(\E)$ of all
endomorphisms of $\E$ is naturally isomorphic to the tensor product
\begin{equation}\label{E:End}
        \End(\E) \ \cong \ C(M)\otimes \EndC(\E).
\end{equation}

By Proposition~3.43 of \cite{BeGeVe}, $F^\E$ decomposes with respect to
\refe{End} as
\begin{equation}\label{E:curv}
        F^\E \ = \ R^\E \ + \ \fes, \qquad\qquad R^\E\in \A^2(M,C(M)), \
                                                \fes\in \A^2(M,\EndC(\E)).
\end{equation}
In this formula, $\fes$ is an invariant of $\n^\E$ called the {\em twisting
curvature \/} of $\E$, and $R^\E$ is determined by the Riemannian curvature $R$
of
 $M$. If  $(e_1\nek e_{2n})$  is an orthonormal frame of the tangent space
$T_xM, \ x\in M$ and $(e^1\nek e^{2n})$ is the dual frame of the cotangent space
 $T^*M$, then
\[
        R^\E(e_i,e_j) \ = \
          \frac14\, \sum_{k,l}\, \langle R(e_i,e_j)e_k,e_l\rangle \, c(e^k)\, c(e^l).
\]

Assume that $\calS$ is a spinor bundle, $\E=\calW\otimes\calS$ and the
connection $\n^\E$ is given by \refe{nE=}. Then $\A(M,\EndC(\E))\cong
\A(M,\End(\calW))$. The twisting curvature   $\fes$   is equal to the curvature
$F^\calW=(\n^\calW)^2$ via this isomorphism (cf. \cite[p.~121]{BeGeVe}). This
explains why $\fes$ is called the twisting curvature.

Let $\calW$ be a vector bundle over $M$ with connection $\n^\calW$ and let
$F^\calW= (\n^\calW)^2$ denote the curvature of this connection. The twisting
curvature of the connection \refe{nEW} on the tensor product $\E\otimes\calW$ is
the sum
\begin{equation}\label{E:WE/S}
        F^{(\E\otimes\calW)/\calS} \ = \ F^\calW \ + \ F^{\E/\calS}.
\end{equation}

\ssec{lichn}{The Lichnerowicz formula}
Let $\E$ be a Clifford module endowed with a Hermitian structure and let
$D:\gme\to \gme$ be a self-adjoint Dirac operator associated to a Hermitian
Clifford connection $\n^\E$. Consider the metric Laplacian (cf. \refss{lapl})
\[
        \Del^\E \ = \ (\n^\E)^*\, \n^\E\, : \ \gme \ \to \gme,
\]
where $(\n^\E)^*$ denotes the operator adjoint to $\n^\E:\gme\to
\Gam(M,T^*M\otimes\E)$   with respect to the  $L_2$-scalar product. Clearly
$\Del^\E$ is a non-negative self-adjoint operator.

The following {\em Lichnerowicz formula \/} (cf. \cite[Theorem~3.52]{BeGeVe})
plays a crucial role in our proof of vanishing theorems:
\begin{equation}\label{E:lichn}
        D^2 \ = \ \Del^\E \ + \ \bfc(\fes) \ + \ \frac{r_M}4,
\end{equation}
where $r_M$ stands for the scalar curvature of $M$ and $\fes$ is the twisting
curvature of $\n^\E$, cf. \refss{curv}. The operator $\bfc(\fes)$ is defied in
\refe{bfc} (see also \refe{bfc2}).

Let $\L$ be a Hermitian line bundle over $M$ endowed with a Hermitian connection
$\n^\L$ and let $\n_k= \n^{\E\otimes\Lk}$ denote the product connection (cf.
\refe{nEW}) on the tensor product $\E\otimes\Lk$. It is a Hermitian Clifford
connection on $\E\otimes\Lk$. We denote by $D_k$ and $\Del_k$ the Dirac operator
and the metric Laplacian associated to this connection. By \refe{WE/S}, it
follows from the Lichnerowicz formula \refe{lichn}, that
\begin{equation}\label{E:lechnW}
        D_k^2 \ = \ \Del_k \ + \ k\, \bfc(F^\L) \ + \ A,
\end{equation}
where $F^\L=(\n^\L)^2$ is the curvature of $\n^\L$ and
\begin{equation}\label{E:A}
        A \ := \ \bfc(\fes) \ + \ \frac{r_M}4 \ \in \ \End(\E) \ \subset \
                \End(\E\otimes\Lk)
\end{equation}
is independent of $\L$ and $k$.

\ssec{cFL}{Calculation of  $\bfc(F^\L)$}
To compare $D_k^2$ with the Laplacian $\Del_k$ we now need to calculate the
operator $\bfc(F^\L)\in \End(\E)\subset \End(\E\otimes\Lk)$. This may be
reformulated as the following problem of linear algebra.

Let $V$ be an oriented Euclidean vector space of real dimension $2n$ and let
$V^*$ denote the dual vector space. We denote by $C(V)$ the Clifford algebra of
$V^*$. Let $E$ be a module over $C(V)$. We will assume that $E$ is endowed with
a Hermitian scalar product such that the operator $c(v):E\to E$ is
skew-symmetric for any $v\in V^*$. In this case we say that $E$ is a {\em
self-adjoint \/} Clifford module over $V$.

The space $E$ possesses a {\em natural grading \/} $E=E^+\oplus E^-$, where
$E^+$ and $E^-$ are the eigenspaces of the chirality operator with eigenvalues
$+1$ and $-1$ respectively, cf. \refss{chir}.

In our applications $V$ is the tangent space $T_xM$ to $M$ at a point $x\in M$
 and  $E$  is the fiber of  $\E$  over  $x$.

Let $F$ be an imaginary valued antisymmetric bilinear form on $V$. Then $F$ may
be considered as an element of $V^*\wedge V^*$. We need to estimate the operator
 $\bfc(F)\in\End(E)$. Here  $\bfc:\Lam V^*\to C(V)$  is the quantization
map defined exactly as in \refss{quant} (cf. \cite[\S{3.1}]{BeGeVe}).

Let us define the skew-symmetric linear map $\tJ:V\to V$ by the formula
\[
        iF(v,w) \ = \ \<v,\tJ w\>, \qquad v,w\in V.
\]
The eigenvalues of $\tJ$ are purely imaginary. Let $\mu_1\ge\cdots\ge\mu_l>0$ be
the positive numbers such that $\pm{i}\mu_1\nek\pm{i}\mu_l$ are all the non-zero
eigenvalues of $\tJ$. Set
\[
        \tau \ = \ \Tr^+ \tJ\ := \ \mu_1+\cdots+\mu_l, \qquad
        m    \ = \ \min_{j}\, \mu_j.
\]
By the Lichnerowicz formula \refe{lichn}, \refp{D-D} is equivalent to the following
\prop{cFL}
Suppose that the bilinear form $F$ is non-degenerate. Then it defines an
orientation of $V$. If this orientation coincides with (resp. is opposite to)
the given orientation of $V$, then the restriction of $\bfc(F)$ onto $E^-$
(resp. $E^+$) is greater than $-(\tau-2m)$, i.e., for any $\alp\in E^-$ (resp.
$\alp\in E^+)$
\[
        \<c(F)\alp,\alp\> \ \ge \ -(\tau-2m)\, \|\alp\|^2.
\]
\eprop
We will prove the proposition in \refss{prcFL} after introducing some additional
constructions. Since we need these constructions also for the proof of
\refp{D-Dac}, we do not assume that $F$ is non-degenerate unless this is stated
explicitly.

\ssec{cs}{A choice of a complex structure on  $V$}
By the Darboux theorem (cf. \cite[Theorem~1.3.2]{Audin91}), one can choose an
orthonormal basis $f^1\nek f^{2n}$ of $V^*$, which defines the positive
orientation of $V$ (i.e., $f^1\wedge\dots\wedge f^{2n}$ is a positive volume form
on $V$) and such that
\begin{equation}\label{E:iF}
       iF_x^\L \ = \  \sum_{j=1}^l\, r_j\, f^j\wedge f^{j+n},
\end{equation}
for some integer $l\le n$ and some non-zero real numbers $r_j$. We can and we
will assume that $|r_1|\ge|r_2|\ge\cdots\ge|r_l|$.

Let $f_1\nek f_{2n}$ denote the dual basis of $V$.
\rem{cFL}
If the vector space $V$ is endowed with a complex structure $J:V\to V$
compatible with the metric (i.e., $J^*=-J$) and such that $F$ is a $(1,1)$ form
with respect to $J$, then the basis $f_1\nek f_{2n}$ can be chosen so that
$f_{j+n}=Jf_j, \ i=1\nek n$.
\erem

Let us define a complex structure $J:V\to V$ on $V$ by the condition
$f_{i+n}=Jf_i, \ i=1\nek n$. Then, the complexification of $V$ splits into the
sum of its holomorphic and anti-holomorphic parts
\[
        V\otimes\CC \ = \ V^{1,0}\oplus V^{0,1},
\]
on which $J$ acts by multiplication by $i$ and $-i$ respectively. The space
$V^{1,0}$ is spanned by the vectors $e_j=f_j-if_{j+n}$, and the space $V^{0,1}$
 is spanned by the vectors  $\oe_j=f_j+if_{j+n}$. Let  $e^1\nek e^n$
and $\oe^1\nek \oe^n$ be the corresponding dual base of $(V\ha)^*$ and
$(V\ah)^*$ respectively. Then \refe{iF} may be rewritten as
\[
        iF_x^\L \ = \  \frac{i}2\, \sum_{j=1}^n\, r_j\, e^j\wedge \oe^j.
\]
We will need the following simple
\lem{mur}
Let $\mu_1\nek\mu_l$ and $r_1\nek r_l$ be as above. Then $\mu_i=|r_i|$, for any
$i=1\nek l$. In particular,
\[
        \Tr^+\tJ \ = \ |r_1| +\cdots+|r_l|.
\]
\elem
\prf
Clearly, the vectors $e_1\nek e_n; \oe_1\nek\oe_n$ form a basis of eigenvectors
of $\tJ$ and
\begin{align}
        \tJ\, e_j \ &= \ ir_j\, e_j, \quad  \tJ\, \oe_j \ = \ -ir_j\, \oe_j \qquad
                        &\mbox{for} \quad j&=1\nek l, \notag\\
        \tJ\, e_j \ &= \ \tJ\, \oe_j \ = \ 0
                        \qquad &\mbox{for} \quad j&=l+1\nek n.\notag
\end{align}
Hence, all the nonzero eigenvalues of $\tJ$ are $\pm{i}|r_1|\nek\pm{i}|r_l|$.
\eprf

\ssec{spin1}{Spinors}
Set
\begin{equation}\label{E:spin1}
        S^+ \ = \ \bigoplus_{j \ \text{even}}\Lam^j(V\ah), \quad
         S^- \ = \ \bigoplus_{j \ \text{odd}}\Lam^j(V\ah).
\end{equation}
Define a graded action of the Clifford algebra $C(V)$ on the graded space
$S=S^+\oplus S^-$ as follows (cf. \refss{almcomp}): \ if $v\in V$ decomposes as
 $v=v\ha+v\ah$  with  $v\ha\in V\ha$  and  $v\ah\in V\ah$, then its
Clifford action on $\alp\in E$ equals
\begin{equation}\label{E:clact}
        c(v)\alp \ = \ \sqrt{2}\, \left( v\ah\wedge\alp  \ - \ \iot(v\ha)\, \alp\right).
\end{equation}
Then (cf. \cite[\S3.2]{BeGeVe}) $S$ is the {\em spinor representation \/} of
$C(V)$, i.e., the complexification $C(V)\otimes\CC$ of $C(V)$ is isomorphic to
$\End(S)$. In particular, the Clifford module $E$ can be decomposed as
\[
        E \ = \ S\otimes W,
\]
where $W=\Hom_{C(V)}\, (S,E)$. The action of $C(V)$ on $E$ is equal to $a\mapsto
c(a)\otimes1$, where $c(a) \ (a\in C(V))$ denotes the action of $C(V)$
on $S$. The natural grading on $E$ is given by $E^\pm= S^\pm\otimes{W}$.

To prove \refp{cFL} it suffices now to study the action of $\bfc(F)$ on $S$. The
latter action is completely described be the following
\lem{ecF}
The vectors $\oe^{j_1}\wedge\dots\wedge\oe^{j_m}\in S$ form a basis of
eigenvectors of $\bfc(F)$ and
\[
        \bfc(F)\, \oe^{j_1}\wedge\dots\wedge\oe^{j_m} \ = \
                \big(\sum_{j'\not\in\{j_1\nek j_m\}} \, r_{j'} \ - \
                        \sum_{j''\in\{j_1\nek j_m\}} \, r_{j''}
                \big)\, \oe^{j_1}\wedge\dots\wedge\oe^{j_m}.
\]
\elem
\prf
Obvious.
\eprf

\ssec{prcFL}{Proof of \refp{cFL}}
Recall that the orientation of $V$ is fixed and that we have chosen the basis
$f_1\nek f_{2n}$ of $V$ which defines the same orientation. Suppose now that the
bilinear form $F$ is non-degenerate. Then $l=n$ in \refe{iF}. It is clear, that
the orientation defined by $iF$ coincides with the given orientation of $V$ if
and only if the number $q$ of negative numbers among $r_1\nek r_n$ is even.
Hence, by \refl{ecF}, the restriction of $\bfc(F)\in\End(S)$ on
$\Lam^j(V\ah)\subset S$ is greater than $-(\tau-2m)$ if the parity of $j$ and
$q$ are different. The \refp{cFL} follows now from \refe{spin1}.
\hfill$\square$

\ssec{prD-Dac}{Proof of \refp{D-Dac}}
Assume that at least $q$ of the numbers $r_1\nek r_l$ are negative and let
$m_q>0$ be the minimal positive number such that at least $q$ of these numbers
are not greater than $-m_q$. It follows from \refl{ecF}, that
\[
        \<\, c(F)\alp,\alp\, \> \ \ge \ -\, (\tau-2m_q)\, \|\alp\|^2,
\]
for any $j<q$ and any $\alp\in \Lam^j(V\ah)$. \refp{D-Dac} follows now
from the Lichnerowicz formula \refe{lichn}. \hfill $\square$

\ssec{prgr}{Proof of \refp{gr}}
Let $\pi_i:\A^{0,*}(M,\E\otimes\Lk)\to \A^{0,i}(M,\E\otimes\Lk)$ denote
the projection and set 
\[
        \tilnk \ = \ \sum_i\, \pi_i\circ\nk\circ\pi_i.
\]
Denote $\tildk=(\tilnk)^*\tilnk$. 
Clearly, $\tildk$ preserves the $\ZZ$-grading on $\A^{0,*}(M,\E\otimes\Lk)$.
It follows from the proof of Theorem~2.16 in
\cite{Dem85}, that there exists a sequence $\eps_1,\eps_2,\ldots$, convergent to
zero, such that
\[
        (1-\eps_k)\, \<\, \tildk\, \gam,\, \gam\, \> \ - \ \eps_kk\, \|\gam\|^2
        \ \le \
        \<\, \Del_k\, \gam,\,  \gam\, \>
        \ \le \
        (1+\eps_k)\, \<\, \tildk\, \gam,\, \gam\, \>
                \ + \ \eps_kk\, \|\gam\|^2,
\]
for any $\gam$ in the domain of $\Del_k$. Hence,
\begin{multline}\notag
        \|\<\, (\Del_k-\tildk)\, \gam,\,  \gam\, \>\| \\
        \ \le \
        \eps_k\, \<\, \tildk\, \gam,\, \gam\, \> \ + \ \eps_kk\, \|\gam\|^2
        \ \le \
        \eps_k \, \sum_{i}\, \<\, \Del_k\, \pi_i\gam,\,  \pi_i\gam\, \>
         \ + \ \eps_kk\, \|\gam\|^2.
\end{multline}

Suppose now that $\gam=\alp+\bet$, where
$\alp\in\Dom^i(D^2_k),\bet\in\Dom^j(D^2_k), \ i\not=j$. Then
\begin{multline}\notag
         \|\, 2\RE\<\, \Del_k\, \alp,\,  \bet\, \>\, \|
        \ = \
         \|\, 2\RE\<\, (\Del_k-\tildk)\, \alp,\,  \bet\, \>\, \| \\
        \ \le \
         \|\, \<\, (\Del_k-\tildk)\, \gam,\,  \gam\, \>\, \|
        \ + \
         \|\, \<\, (\Del_k-\tildk)\, \alp,\,  \alp\, \>\, \|
        \ + \
         \|\, \<\, (\Del_k-\tildk)\, \bet,\,  \bet\, \>\, \| \\
        \ \le \
         2\eps_k\, \<\, \Del_k\, \alp,\,  \alp\, \>
        \ + \
         2\eps_k\, \<\, \Del_k\, \bet,\,  \bet\, \>
        \ + \
         2\eps_kk\, \|\alp\|^2 \ + \ 2\eps_kk\, \|\bet\|^2.
\end{multline}
Similarly one obtains an estimate for the imaginary part of
$\<\Del_k\alp,\bet\>$. This leads to the following analogue of \refp{gr} for the
operator $\Del_k$:
\begin{equation}\label{E:albe}
         \|\, \<\, \Del_k\, \alp,\,  \bet\, \>\, \|
        \ \le  \
         2\eps_k\, \<\, \Del_k\, \alp,\,  \alp\, \>
        \ + \
         2\eps_k\, \<\, \Del_k\, \bet,\,  \bet\, \>
        \ + \
         2\eps_kk\, \|\alp\|^2 \ + \ 2\eps_kk\, \|\bet\|^2.
\end{equation}

We now apply the Lichnerowicz formula \refe{lechnW} to obtain \refp{gr} from
\refe{albe}. Note, first, that the operator
$A\in\End(\E)\subset\End(\E\otimes\Lk)$, defined in \refe{A}, is independent of
$k$ and bounded. Note, also, that, by \refl{ecF}, the operator $\bfc(F^\L)$
preserves the $\ZZ$-grading on $\A^{0,*}(M,\Ek)$. Hence, it follows from
\refe{albe} and the Lichnerowicz formula \refe{lechnW} that
\begin{multline}\notag
          \|\, \<\, D_k^2\, \alp,\,  \bet\, \>\, \|
         \ \le \
          |\, \<\, \Del_k\, \alp,\,  \bet\, \>\, |
          \ + \
                   |\, \<\, A\, \alp,\,  \bet\, \>\, | \\
         \ \le \
         2\eps_k\, \<\, \Del_k\, \alp,\,  \alp\, \>
         \ + \
         2\eps_k\, \<\, \Del_k\, \bet,\,  \bet\, \>
         \ + \
         2\eps_kk\, \|\alp\|^2 \ + \ 2\eps_kk\, \|\bet\|^2
         \ + \
         \|A\|\, \|\alp\|\, \|\bet\| \\
        \ \le \
         2\eps_k\, \<\, D_k^2\, \alp,\,  \alp\, \>
         \ + \
         2\eps_k\, \<\, D_k^2\, \bet,\,  \bet\, \>
         \ + \
         2\eps_kk\, \big(1+\|\bfc(F^\L)\|+2\|A\|\big)\, \|\alp\|^2 \\
         \ + \
         2\eps_kk\, \big(1+\|\bfc(F^\L)\|+2\|A\|\big)\, \|\bet\|^2.
\end{multline}
Hence, \refp{gr} holds with $\del_k=(1+\|\bfc(F^\L)\|+2\|A\|\big)\eps_k$.
\hfill$\square$

\sec{lapl}{Estimate of the metric Laplacian}

In this section we prove \refp{lapl}.
\ssec{reduc}{Reduction to a scalar operator}
In this subsection  we construct a space $\Z$
and an operator $\tilDel$ on the space $L_2(\Z)$ of $\Z$, such that the operator
$\Del_k$ is ``equivalent" to a restriction of $\tilDel$ onto certain subspace of
$L_2(\Z)$. This allow to compare the operators $\Del_k$ for different
values of $k$.

Let $\F$ be the principal $G$-bundle with a compact structure group $G$,
associated to the vector bundle $\E\to M$. Let $\Z$ be the principal
$(S^1\times{G})$-bundle over $M$, associated to the bundle $\E\otimes\L\to M$.
Then $\Z$ is a principle $S^1$-bundle over $\F$. We denote by $p:\Z\to \F$ the
projection.

The connection $\n^\L$ on $\L$ induces a connection on the bundle
$p:\Z\to\F$. Hence, any vector $X\in T\Z$ decomposes as a sum
\begin{equation}\label{E:dec}
        X \ = \ X{\hor} \ + \ X\vert,
\end{equation}
of its horizontal and vertical components.

Consider the {\em horizontal exterior derivative \/} $d\hor:
C^\infty(\Z)\to \A^1(\Z,\CC)$, defined by the formula
\[
        d\hor f(X) \ = \ df(X\hor), \qquad X\in T\Z.
\]

The connections on $\E$ and $\L$, the Riemannian metric on $M$, and
the Hermitian metrics on $\E, \L$ determine a natural Riemannian
metrics $g^\F$ and $g^\Z$ on $\F$ and $\Z$ respectively,
cf. \cite[Proof of Theorem~2.1]{BoUr96}. Let $(d\hor)^*$ denote the
adjoint of $d\hor$ with respect to the scalar products induced by this
metric. Let
\[
        \tD \ = \ (d\hor)^*d\hor:\, C^\infty(\Z) \to C^\infty(\Z)
\]
be the {\em horizontal Laplacian} for the bundle $p:\Z\to \F$.

Let $C^\infty(\Z)_k$ denote the space of smooth functions on $\Z$, which are
homogeneous of degree $k$ with respect to the natural fiberwise circle action on
the circle bundle $p:\Z\to \F$. It is shown in \cite[Proof of
Theorem~2.1]{BoUr96}, that to prove \refp{lapl} it suffices to prove \refe{lapl}
for the restriction of $\tD$ to the space $C^\infty(\Z)_k$.

\ssec{symbol}{The symbol of $\tD$}
The decomposition \refe{dec} defines a splitting of the cotangent bundle $T^*\Z$
to $\Z$ into the horizontal and vertical subbundles. For any $\xi\in T^*\Z$, we
denote by $\xi\hor$ the horizontal component of $\xi$. Then, one easily checks
(cf. \cite[Proof of Theorem~2.1]{BoUr96}), that the principal symbol
$\sig_2(\tD)$ of $\tD$ may be written as
\begin{equation}\label{E:symb}
        \sig_2(\tD)(z,\xi) \ = \ g^\F(\xi\hor,\xi\hor).
\end{equation}
The subprincipal symbol of $\tD$ is equal to zero.

On the {\em character set \/} $\calC=\big\{(z,\xi)\in T^*\Z\backslash\{0\}:\,
\xi\hor=0\big\}$ the principal symbol $\sig_2(\tD)$ vanishes to second order.
Hence, at any point $(z,\xi)\in\calC$, we can define the {\em Hamiltonian map
\/} $F_{z,\xi}$ of $\sig_2(\tD)$, cf. \cite[\S21.5]{Horm3}. It is a
skew-symmetric endomorphism of the tangent space $T_{z,\xi}(T^*\Z)$. Set
\[
        \Tr^+ F_{z,\xi} \ = \ \nu_1 \ + \ \cdots \ + \ \nu_l,
\]
where $i\nu_1\nek i\nu_l$ are the nonzero eigenvalues of $F_{z,\xi}$
for which $\nu_i>0$.

Let $\rho:\Z\to M$ denote the projection. Then, cf. \cite[Proof of
Theorem~2.1]{BoUr96}
\footnote{The absolute value sign of $\xi\vert$ is erroneously missing
in \cite{BoUr96}.},
\begin{equation}\label{E:tr+F}
        \Tr^+ F_{z,\xi} \ = \ \tau(\rho(z))\, |\xi\vert|
\end{equation}
Here   $\xi\vert$ is the vertical component of $\xi\in T^*\Z$, and
$\tau$ denotes the function defined in \refe{tau}.

\ssec{melin}{Application of the Melin  inequality}
Let $D\vert$ denote the generator of the $S^1$ action on $\Z$. The
symbol of $D\vert$ is $\sig(D\vert)(z,\xi) =\xi\vert$. Fix $\eps>0$,
and consider the operator
\[
        A \ = \ \tD \ - \big(\tau(\rho(z))-\eps\big)\, D\vert:\,
                C^\infty(\Z) \ \to \  C^\infty(\Z).
\]
The principal symbol of $A$ is given by \refe{symb}, and the subprincipal symbol
\[
        \sig_1^s(A)(z,\xi) \ = \  -\big(\tau(\rho(z))-\eps\big)\, \xi\vert.
\]
It follows from \refe{tr+F}, that
\[
        \Tr^+ F_{z,\xi}+\sig_1^s(A)(z,\xi) \ \ge \
                \eps\, |\xi\vert| \ > \ 0.
\]
Hence, by the Melin inequality (\cite{Melin71}, \cite[Theorem~22.3.3]{Horm3}),
there exists a constant $C_\eps$, depending on $\eps$, such that
\begin{equation}\label{E:Aff}
        \<\, Af,f\, \> \ \ge -\, C_\eps\, \|f\|^2.
\end{equation}
Here $\|\cdot\|$ denotes the $L_2$ norm of the function $f\in
C^\infty(\Z)$.

{}From \refe{Aff}, we obtain
\[
        \<\, \tD f,f\, \> \ \ge \
                \<\, (\tau(\rho(z))-\eps)D\vert f,f\, \> \
        - \  C_\eps\, \|f\|^2.
\]
Noting that if $f\in  C^\infty(\Z)_k$, then $D\vert f=kf$, the proof
is complete. \hfill $\square$

\providecommand{\bysame}{\leavevmode\hbox to3em{\hrulefill}\thinspace}

\end{document}